\documentclass{amsart}
\def\tensor{\,\raise2pt\hbox{${}_{\otimes}$}\,}% Tensor
\def\fdg{\,:\,}% ... fuer die gilt ...
\def\ptl{\partial}% Partial
\def\rest#1{\raise-2pt\hbox{${\lfloor_{#1}}$}}% Restringiert
\def\ip#1#2{\langle#1,#2\rangle}
\def\olin#1{\overline{#1}{}}% Oben-quer
\def\tild#1{\widetilde{#1}{}}% Tilde
\def\grad{{\nabla}}% Gradient
\newcommand{\leftexp}[2]{{\vphantom{#2}}^{#1}{#2}}
\def\halb{\frac{1}{2}}% 1/2
\def \a{\alpha}

\newtheorem{theorem}{Theorem}[section]
\newtheorem{lemma}[theorem]{Lemma}

\newtheorem{remark}[theorem]{Remark}
\newtheorem{definition}[theorem]{Definition}
\newcommand{\ba}{\begin{array}}
\newcommand{\ea}{\end{array}}
\newcommand{\bea}{\begin{eqnarray}}
\newcommand{\eea}{\end{eqnarray}}
\newcommand{\bee}{\begin{eqnarray*}}
\newcommand{\eee}{\end{eqnarray*}}

\renewcommand{\a}{\alpha}

\usepackage{color}

\newcommand{\green}[1]{{\color{green}#1}}

\newcounter{mnotecount}[section]

\renewcommand{\themnotecount}{\thesection.\arabic{mnotecount}}

\newcounter{mymnotecount}[section]
\renewcommand{\themymnotecount}{\thesection.\arabic{mymnotecount}}
\newcommand{\mymnote}[1]{\protect{\stepcounter{mymnotecount}}${\raisebox{0.5\baselineskip}[0pt]{\makebox[0pt][c]{\color{green}{\tiny\em$\bullet$\themnotecount}}}}$\marginpar{\raggedright\tiny\em$\!\bullet$\themymnotecount:

\green{#1}}\ignorespaces}

\renewcommand{\mymnote}[1]{}
\begin{document}

\title{ A Note on The Dimensional Reduction of Axisymmetric spacetimes}


%    author one information
\author{Nishanth Gudapati}
\address{Department of Mathematics, Yale University, 10 Hillhouse Avenue, New Haven, CT-06511, USA }
\curraddr{}
\email{nishanth.gudapati@yale.edu}
\thanks{ (NG) is supported by a Deutsche Forschungsgemeinschaft (DFG) Postdoctoral Fellowship GU1513/1-1 to Yale University}

%    \subjclass is required.
\subjclass[2010]{Primary: 83, Secondary 35}

%\date{\today}

%    Abstract is required.
\iffalse
\begin{abstract}
The aim of this notes is to associate the ADM mass of axisymmetric spacetimes with the 2+1 Ashtekar-Varadarajan mass via the usual dimensional reduction procedure (ala Kaluza-Klein). In this connection, a by product is to understand the induced spatial asymptotic flatness together with the matter fields and rectify the degeneracies of the energy in the 2+1 dimensional scenario. Separately, we also aim to understand the positivity of `geoemtric-mass' for Einstein's equations with $\Lambda >0$, to assist in the larger program initiated by Ashtekar et. al. 
\end{abstract}
\fi
\begin{abstract}
 We investigate the dimensional reduction of 3+1 vacuum axisymmetric Einstein's equations to 2+1 dimensional Einstein-wave map system and observe that the resulting system is 1) not asymptotically flat, 2) its geometric-mass diverges and 3) the energy of wave map also diverges. Subsequently, we discuss the consequences of these issues.
\end{abstract}
\maketitle

\section{Background and Preliminaries}
Let $(M, g)$ be a 2+1 dimensional globally hyperbolic, Lorentzian manifold and further suppose that $(M, g)$ admits a 2+1 decomposition such that it is foliated by a family of Cauchy hypersurfaces $(\Sigma_t, q_t)$. Let us assume that each $(\Sigma_t, q_t)$ is asymptotically flat in the sense of Ashtekar-Varadarajan \cite{ash_var}:
\begin{align}
q_{ab} = r^{-m_{AV}} (\delta_{ab} + \mathcal{O} (r^{-1})),
\end{align}   
in the asymptotic region where $\delta_{ab}$ is the 2d flat metric.
%\mnote{check!}
%\[ \delta_{ij} = D_i x D_j x + D_i y D_j y.\]

\noindent It was long known that the notion of energy in general relativity is most reasonable for an isolated self-gravitating system that converges to the trivial solution asymptotically. In 1960s, following contributions from Dirac, ADM have demonstrated that energy at an instance can be measured in terms of the deviation from this trivial solution in the asymptotic region. The underlying notion is that the curvature of the asymptotically flat 3-manifold introduces a deviation from flat space which can be measured canonically at infinity with the trivial solution as the reference.  Several fundamental results have followed which established positivity, consistency and convergence of this energy\cite{schoen-yau-1, schoen-yau-2,witten-pmt,bartnik-adm}.  In another foundational result, which shall play a crucial role in our analysis, Ashtekar and Varadarajan \cite{ash_var} have established that $m_{AV}$ is the natural Hamiltonian analogue of ADM-mass for the 2+1 case and that $m_{AV}$ is non-negative with $m_{AV} =0$ iff $q$ is flat. %The work of Ashtekar and Varadarajan establishes that, for an asymptotically flat 2-manifold, the geometric mass ($m_{AV}$) in the asymptotic region can be measured in terms of the angle deficit at infinity. 

\noindent Now suppose we are interested in the Einstein-wave map system in 2+1 dimensions. As we shall see later,  the reason we are considering 2+1 Einstein-wave system is that it occurs in vacuum 3+1 Einstein's equations.
\begin{align}
\mathbf{E}_{\mu \nu} =& \mathbf{T}_{\mu \nu} \quad \textnormal{on} \quad (M, g)  \notag \\
\square_g U^i + \leftexp{(h)}{\Gamma}^i_{jk} (U) g^{\mu \nu}\ptl_\mu U^j \ptl_\nu U^k =&0 \quad \quad\,\, \textnormal{on} \quad (M, g),
\end{align}
where $U$ is the wave map 
\begin{align}
U \fdg (M, g) \to (N, h),
\end{align}

$\mathbf{E}$ is the  Einstein-tensor of $(M, g)$, $\mathbf{T}$ is the energy-momentum tensor for wave maps:
\[ \mathbf{T}_{\mu \nu} = \ip{\ptl_\mu U}{ \ptl_\nu U}_h - \halb g_{\mu \nu} \ip{\ptl^\sigma U}{\ptl_\sigma U}    \]
$\square_g$ is the covariant wave operator, $\Gamma$'s are the Christoffel symbols of the target $(N, h)$. For demonstration purposes, let us restrict to the case where $(N, h)$ is a hyperbolic 2-plane, which can also be interpreted a surface of revolution with the generating function $f$. As in \cite{diss}, let us further restrict to the case where $(M, g)$ is rotationally symmetric and the map $U$ is an equivariant wave map, which roughly means that the map commutes with the rotational symmetries of $(M, g)$, $ (N, h)$ and thus partially decouples the wave maps equation to give the system:
\begin{align}\label{equi-ewm}
\mathbf{E}_{\mu \nu} =& \mathbf{T}_{\mu \nu}  \notag \\
\square_g u =& \frac{f_u (u) f(u)}{r^2}.
\end{align}

Suppose $(M, g)$ can be expressed in the following form in polar coordinates $(t, r, \theta)$:
\begin{align}
ds^2_g = -e^{2\Omega (t, r)} dt^2 + e^{2\gamma(t, r)} dr^2 + r^2 d\theta^2
\end{align}
with $\gamma, \ptl_r \gamma$ and $\Omega =0$ on this axis $\Gamma.$
The spatial part is 
\begin{align}
q = e^{2 \gamma} dr^2 + r^2 d\theta^2 
\end{align}
with $r \in [0, \infty)$ and $\theta \in [0, 2\pi).$
Suppose 
\[ \gamma_{\infty} = \lim_{r \to \infty} \gamma \]
In view of our assumption on asymptotic flatness of $q$ the limit above converges.
Therefore, the form of the metric at the spatial asymptotic region is given by 
\begin{align}
q_{\infty} \fdg = e^{2 \gamma_{\infty}} dr_{\infty}^2 + r_{\infty}^2 d\theta^2  
\end{align}
where $r_\infty$ is the 2d Euclidean radius function.
\subsection*{Angle-deficit}
In order to compare the asymptotic behaviour of our metric at spatial infinity to the flat metric let us 
introduce the variables
\[ \rho \fdg = e^{\gamma_\infty} r_\infty\quad  \textnormal{and}\quad  \vartheta \fdg = e^{-\gamma_\infty} \theta  \]
with 
\begin{align}
0\leq \vartheta < 2\pi e^{-\gamma_\infty}.
\end{align}
So that $q_{\infty}$ can be represented in the form of the flat metric in the asymptotic region as
\begin{align}
q_{\infty} = d\rho^2 + \rho^2 d\vartheta^2, \quad 0\leq \vartheta < 2\pi e^{-\gamma_\infty}.
\end{align}
Thus the angle deficit at infinity is given by 
\begin{align}
\ptl \theta_{\infty} \fdg = 2\pi (1-e^{-\gamma_\infty}) 
\end{align}
As we shall see next, this notion is streamlined by Ashtekar-Varadarajan \cite{ash_var}.
%We shall see later that $\gamma$, which is $0$ on the axis, is monotonically increasing function of $r.$
\subsection*{AV-Mass}
In order the deduce the AV-mass of our system, we need to rewrite our form of the metric 
in the AV-form. Let us define the asymptotic Cartesian coordinates $x = r_\infty \cos \theta$ and $y =r_\infty \sin \theta$, so that the flat metric is 
\begin{align}
\delta_{ab} = \grad_a x \grad_b x + \grad_a y \grad_b y
\end{align}

\noindent Let us introduce the variable $\varrho$ such that our metric $q$ can be expressed in the form 
\[ q_{ab} = \varrho^{-m_{AV}} (\delta_{ab} + \mathcal{O} (\varrho^{-1}))\]
in the asymptotic region. 
Suppose
\begin{align}
\varrho^{2-m_{AV}} = r_{\infty}^2
\end{align}  
so that $dr_{\infty} = \left(\frac{2-m_{AV}}{2} \right) \varrho^{-\frac{m_{AV}}{2}} d\varrho$ which implies that $e^{2\gamma_\infty} = \frac{2}{2-m_{AV}},$ so
\begin{align}
m_{AV} = 2(1-e^{-\gamma_{\infty}}).
\end{align}
We shall demonstrate later that $\gamma$, which is 0 on the axis $\Gamma$, is a monotonically increasing function with respect to $r$. This means that the AV-mass is also bounded from above $m_{AV} <2.$
\subsection*{Energy}
Let $\mathcal{N}$ is the unit time-like normal of the $(\Sigma, q) \to (M, g)$ embedding, then the Einstein-Hamiltonian constraint is formally given by
\[R_q + (K^a_{\,a} )^2 - K_{ab} K^{ab} = \mathbf{T} (\mathcal{N}, \mathcal{N}) \]
which is essentially the equation, 
\[  \mathbf{E} (\mathcal{N}, \mathcal{N}) = \mathbf{T} (\mathcal{N}, \mathcal{N}) \]
reduced by the Gauss-Kodazzi equation. The Einstein-Hamiltonian constraint shall play an important role in our analysis as it captures the energy in the system. But before we move on to the energy let us briefly discuss more details about wave maps. Wave maps are critical points of a functional that is a natural geometric generalization of harmonic maps:

\begin{align}
S_{WM} \fdg =& \int_{M} \mathcal{L}_{WM} \bar{\mu}_g \notag\\
\fdg=& \halb \int_{M} g^{\mu \nu} h_{ij} \ptl_\mu U^i \ptl_\nu U^j \bar{\mu}_g.
\end{align}
In the case that we considered in \eqref{equi-ewm}, we have a wave map that is evolving with respect to a time function on the manifold $(M, g)$ which has the wave map itself as the source. As the wave map is propagating in time its canonical stress energies are encloded in the energy-momentum tensor 
\begin{align}
\mathbf{T}_{\mu \nu} \fdg = \frac{\ptl L_{WM}}{\ptl g_{\mu \nu}} - g_{\mu \nu} L_{WM}.
\end{align} Then the energy is given by 
\begin{align}
E(t) \fdg = \int_{\Sigma_t} \mathbf{T} (\mathcal{N}, \mathcal{N})\, \bar{\mu}_q
\end{align}
If we were in the Cartesian Minkowski space $\mathbb{R}^{2+1},$
\begin{align}
E(t) = \halb \int_{\mathbb{R}^2} \Vert \ptl_t U\Vert^2_h + \Vert \grad_x U \Vert^2_h \, d^2 x
\end{align}

In our case, 
\begin{align}
E(t) =& \int_{\Sigma_t} \mathbf{T} (\mathcal{N}, \mathcal{N}) \sqrt{q}\, dr d\theta \notag\\
=& \halb \int_{\Sigma_t} \left( \big\Vert \mathcal{N} U \big\Vert_h + \big\Vert \mathcal{R} U \big\Vert_h + \big\Vert \frac{1}{r}\ptl_\theta U \big\Vert_h \right) \sqrt{q}\, dr d\theta, 
\quad \mathcal{R} \fdg= e^{-\gamma}\ptl_r \notag\\
=& \halb \int_{\Sigma_t} \left( e^{-2 \Omega} (\ptl_t u)^2 + e^{-2 \gamma} (\ptl_r u)^2 + \frac{f^2(u)}{r^2}  \right) \sqrt{q}\, dr d\theta 
\end{align}
Crucially, It should be noted that a priori we can impose conditions on the Cauchy slice so that $E(t)$ converges.  Now let us return to the Hamiltonian constraint: A calculation shows that the Einstein tensor is 
\begin{align}
\mathbf{E} (\mathcal{N}, \mathcal{N}) = e^{-2 \gamma} r^{-1} \ptl_r \gamma
\end{align}
Thus the Hamiltonian constraint gives
\begin{align}
e^{-2 \gamma} r^{-1} \ptl_r \gamma = \halb \left( e^{-2 \Omega} (\ptl_t u)^2 + e^{2 \gamma} (\ptl_r u)^2 + \frac{f^2(u)}{r^2}  \right) 
\end{align}
This means
\begin{align}
e^{-\gamma} \ptl_r \gamma = r e^{ \gamma} \halb \left( e^{-2 \Omega} (\ptl_t u)^2 + e^{2 \gamma} (\ptl_r u)^2 + \frac{f^2(u)}{r^2}  \right) \notag\\
- \ptl_r (e^{-\gamma}) =  r e^{ \gamma} \halb \left( e^{-2 \Omega} (\ptl_t u)^2 + e^{2 \gamma} (\ptl_r u)^2 + \frac{f^2(u)}{r^2}  \right) 
\end{align}
Integrating with respect to $r$, we shall obtain 
\begin{align}\label{gmon}
1-e^\gamma = \frac{1}{4\pi} \int^r_{0} \left( e^{-2 \Omega} (\ptl_t u)^2 + e^{-2 \gamma} (\ptl_r u)^2 + \frac{f^2(u)}{r^2}  \right) \sqrt{q}\, dr d\theta 
\end{align}
The equation shows that $\gamma$ is monotonically increasing with respect to $r$
\[ 0 \leq \gamma \leq \gamma_{
\infty}\]
Further,\eqref{gmon} implies 
\begin{align}
1-e^{\gamma_{\infty}} = \frac{1}{2\pi} E(t)
\end{align}
So that
\begin{align}
e^{\gamma_{\infty} }= \left( 1-\frac{1}{2\pi} E(t) \right)^{-1}
\end{align} 
%The above equation shows that, as $E(t) \to \infty$ the metric function $e^{\gamma_\infty} \to 0$ and $m_{AV} \to 2.$
Combining them all,
\begin{align}
\ptl \theta_\infty = \pi m_{AV} =  E(t).
\end{align}

Thus we have established that there is a neat correlation between all these quantities.  Weaker conservation laws associated to these quantities were useful in \cite{diss} to obtain non-concentration for the system for large data.  %Some aspects of these computations were also done in \cite{laan_equinotes}. %Furthermore, all of the aforementioned quantities can also be represented in terms of an error term in the spatial asymptotic limit of the Gauss-Bonnet theorem \cite{}.  
\\

\noindent In the last few decades there has been an impressive progress in the technical aspects of dispersive and geometric PDE. Given this background, during his PhD work the author had observed and duly reported to the wider community that the rich structure due to the scaling symmetries in the 2+1 Einstein-wave map system can be used to obtain hitherto inaccessible results in the context of the initial value problem of Einstein's equations for general relativity, especially for large data. To this end, the approach was able to produce a large data global result for the equivariant Einstein-wave map system \cite{AGS}, by building on the non-concentration result in \cite{diss}. Such a global result was the first of its kind.  
%To this end, global existence for large data was obtained for equivariant Einstein-map system by building on the preceding nonconcentration result in \cite{diss}, which marked the first instance of a global regularity result for  large data for noncompact Cauchy hypersurfaces with one Killing field.
In theory one can hope to use this formulation
in the dynamical black hole stability problem, particularly since the Kerr family admits a similar dimensional reduction. However, in this note we shall demonstrate that the Einstein-wave map approach to study dynamical axisymmetric spacetimes is fraught with significant foundational issues. These issues are in contrast with the translational symmetry case.  The details are discussed in the next section.

\section*{3+1 Axisymmetric Spacetimes}
Suppose $(\bar{M}, \bar{g})$ is a 3+1 dimensional globally hyperbolic Lorentzian spacetime that admits a foliation of 3d (Riemannian) Cauchy hypersurfaces which are asymptotically flat and its ADM-mass converges:
\iffalse
 following expansion 
 \begin{align}
\bar{q}_{ij} = \left( 1+ \frac{m}{r} \right) \bar{\delta}_{ij} + \mathcal{O} (r^{-1-\a})
\end{align}
where $\bar{\delta}_{ij}$ is the flat space 3-metric expressed in Cartesian coordinates. The quantity $m$ is the ADM-mass of the 3-manifold $\bar{\Sigma}, \bar{q}$, 

which is formally defined as the limit
\fi
\begin{align}\label{defadm}
m_{ADM} \fdg= \frac{1}{16\pi} \lim_{r \to \infty} \int_{\mathbb{S}^2(r)} \bar{\delta}^{kl} \left( \ptl_k \bar{q}_{il} - \ptl_i \bar{q}_{lk}  \right) \frac{x^i}{\bar{r}} \bar{\mu}_{\mathbb{S}^2(r)} 
\end{align}
%It can be proven that, for $\alpha \geq 1/2$ the limit \eqref{defadm} converges. 
Now suppose $(\bar{\Sigma}, \bar{q})$ is spatially rotationally symmetric: $\exists$ a space-like rotational Killing vector that has invariant points ($\Gamma$) and closed integral curves. 

In the presence of Killing vectors, the Einstein's equations admit a reduction to 2+1 Einstein-wave map system. This is also applicable for a rotational Killing vector in the regular region. For demonstration purposes, let us go ahead and perform the dimensional reduction. Consider the following ansatz for $(\bar{M}, \bar{g}):$

\begin{align} \label{axi-ansatz}
\bar{g}_{\mu \nu} dx^\mu dx^\nu = g_{\mu \nu} dx^\mu dx^\nu + e^{2u} \left( d\phi + A_\nu dx^\nu \right)^2 \quad \text{for} \quad \mu, \nu =0,1, 2, 3.
\end{align}
where $\ptl_{\phi}$ is the rotational Killing vector, $g$ is the metric of $(M,g) \fdg (\bar{M}, \bar{g}) \setminus SO(2)$ where the axis is identified by $\Vert \ptl_\phi\Vert_{\bar{g}} =0$ (with the convention $e^{-\infty} \fdg =0.$) In the following $\grad$ denotes the covariant derivative defined in $(M, g).$

With the above ansatz, the 3+1 vacuum Einstein's equations for $(\bar{M}, \bar{g})$ 
\begin{align}
\bar{R}_{\mu \nu} =0
\end{align}  
can be rewritten as 

\begin{subequations} \label{Ricci-high-low}
\begin{align}
0=\bar{R}_{\mu \nu} =& R_{\mu \nu} - \grad_\mu u \grad_\nu u - \grad_{\mu} \grad_ {\nu} u - \halb e^{2 u} F_{\mu \sigma} F^\sigma_\nu \label{Rmunu} \\
0=\bar{R}_{\mu 3} =& -\halb e^{-u} \grad_{\sigma} (e^{3 u} F^\sigma _\mu) \label{Rmu3} \\
0=\bar{R}_{3 3} =& -e^{2u} (g^{\mu \nu} \grad_\mu \grad_\nu u + g^{\mu \nu} \grad_\mu u \grad_\nu u -\frac{1}{4} e^{2u} F_{\mu \nu} F^{\mu \nu} ). \label{R33}
\end{align}
\end{subequations}
where  
\[ F_{\mu \nu} \fdg = \grad_{\mu} A_\nu- \grad_{\nu} A_\mu \]

is the Faraday tensor in the chosen coordinate frame. Before we move on, we shall need the following formulas concerning the conformal transformations of Ricci tensor and wave operator.
Suppose
\[ \tilde{g} \fdg = e^{2\psi} g \]
We have,
\begin{subequations}
\begin{align}
\sqrt{-\tilde{g}} =& e^{3 \psi} \sqrt{-g} \\
\square_{\tilde{g}} u =&  \frac{1}{\sqrt{-\tilde{g}}} \ptl_\nu \left( \sqrt{-\tilde{g}}\, \tilde{g}^{\mu \nu}\ptl_\mu u \right) = e^{-2 \psi} \left( \square_g u + g^{\mu \nu} \ptl_\nu \psi \, \ptl_\mu u \right)  \\
\tilde{R}_{\mu \nu} =& R_{\mu \nu} - g_{\mu \nu} \grad^\sigma \grad_\sigma \psi -  \grad_\mu  \grad_\nu \psi + \grad_\mu \psi \grad_\nu \psi - g_{\mu \nu} \grad^\sigma \psi \grad_\sigma \psi. 
\end{align}
\end{subequations}
$\mu, \nu, \sigma =0, 1, 2.$ %These conformal transformations shall simplify the terms in the right side of \eqref{} if we use the Ricci tensor $\tilde{R}_{\mu \nu}$ of $(\tilde{M}, \tilde{g})$

If we examine the equation \eqref{Rmu3} carefully, we can interpret this as the closure of a 1-form. Define

\[ G \fdg = e^{3\psi} \, \leftexp{*}{F} \]
Then \eqref{Rmu3} for $ \bar{R}_{\mu 3}=0$ implies
\[ d\, G =0. \]
This in turn implies there exists a twist potential $v$ such that
\[ G = d\,v \]

The equation $d\, F =0$ for the twist potential $v$ can be transformed into the wave maps equation in the conformally transformed manifold $(\tilde{M}, \tilde{g})$ with $\psi = u.$  
\begin{align} \label{twist_pot}
 \tild{\grad}_{\mu}(e^{-3u} \tild{g}^{\mu \nu} \ptl_\nu v) =0 
\end{align}
This equation constitutes one of the wave maps equations. Subsequently, a linear combination of the equations \eqref{Rmunu} and \eqref{R33} rewritten in the conformally transformed $(\tilde{M}, \tilde{g})$ result in the Einstein's equations and the other wave maps equation:

\begin{subequations}
\begin{align}
\tilde{\mathbf{E}}_{\mu \nu} =&\, \tilde{\mathbf{T}}_{\mu \nu} \\
\tilde{\grad}^\mu \ptl_\mu u + \halb e^{-4u} g^{\mu \nu} \ptl_\mu v \ptl_\nu v =& \,0 \\
\tilde{\grad}^\mu \ptl_\mu v - 4 \, g^{\mu \nu} \ptl_\mu u \ptl_\nu v =&\,0.
\end{align}
\end{subequations}

To simplify our notation from now on, we shall denote $(\tilde{M}, \tilde{g})$ by $(M, g)$ itself. So,

\begin{subequations} \label{axi-ewm}
\begin{align}
\mathbf{E}_{\mu \nu} =& \mathbf{T}_{\mu \nu} \quad \textnormal{on} \quad (M, g)  \notag \\
\square_g U^i + \leftexp{(h)}{\Gamma}^i_{jk} (U) g^{\mu \nu}\ptl_\mu U^j \ptl_\nu U^k =&0 \quad \quad\,\, \textnormal{on} \quad (M, g), \label{axi-wm}
\end{align}
\end{subequations}
where $U$ is now the specific wave map 
\begin{align}
U \fdg (M, g) \to (\mathbb{H}^2, h).
\end{align}

%All of this derivation can be performed more elegantly in variational formulation, but this Euler-Lagrangian version shall me more useful for our discussion. 
Therefore, using the rotational Killing vector $\ptl_{\phi}$ we have superficially arrived at the elegant 2+1 Einstein-wave map system. In 1960s Ernst had arrived at a similar system but with a crucial difference (more later). Now let us perform the reduction for some well-known spacetimes.

\subsection*{Dimensional Reduction of Minkowski}

Consider the 3+1 dimensional Minkowski spacetime:

\begin{align}
ds^2 = -dt^2 + dr^2 + r^2 d\theta^2 + r^2 \sin^2 \theta d\phi^2,
\end{align}
$\theta \in [0, \pi], \phi \in [0, 2\pi)$.
The spatial part is 
\begin{align}
\bar{q} = dr^2 + r^2 d\theta^2 + r^2 \sin^2 \theta d\phi^2
\end{align}
This can be rewritten in the dimensional reduction ansatz:
\begin{align}
\bar{q} = e^{-2u} q + e^{2u} \Phi^2 
\end{align}
where $e^{2u} \fdg= r^2 \sin^2 \theta$. To illustrate our main points, let us cut out  the axis i.e., $u = \log(r\sin \theta)$ and $\Phi$ in hypersurface othogonal case is $d \phi$
\[ q = r^2 \sin^2 \theta (dr^2 + r^2 d\theta^2) \]
The first immediate  observation we have arrived at a manifold that has curvature from the flat Minkowski space using the dimensional reduction procedure. 
Furthermore, we note that this manifold is not asymptotically flat and has infinite curvature at statial infinity. To prove that its geometric-mass also diverges, consider the following equations resulting from the dimensional reduction 
\begin{align}
\mathbf{E}_{\mu \nu} =& \mathbf{T}_{\mu \nu} \quad \textnormal{on} \quad (M, g) \notag \\
\square_g u =&0 \quad \quad \,\, \textnormal{on} \quad (M, g)
\end{align}
where 
\[ g = r^2 \sin^2 \theta ( -dt^2 + dr^2 + r^2 d\theta^2) \]
and
\[\mathbf{T}_{\mu \nu} = \grad_\mu u \grad_\nu u - \halb g_{\mu \nu} \grad^\sigma u \grad_\sigma u.\]

The unit normal of the $\Sigma \hookrightarrow M$ embedding is $\mathcal{N} = (r\sin\theta)^{-1}\ptl_t$ and $\sqrt{q} = r^3 \sin^2 \theta.$ Also note that 

\begin{align}
 \ptl^\sigma u \ptl_\sigma u =& g^{rr}\ptl_r u \ptl_r u + g^{\theta \theta} \ptl_\theta u \ptl_\theta u  \notag \\
 =& \frac{1}{r^4 \sin^4 \theta}
\intertext{and}
\mathbf{T}_{tt} =& - \halb g_{tt} \ptl^\sigma u \ptl_\sigma u \notag \\
=& \frac{1}{2 r^2 \sin^2 \theta}\\
\intertext{thus}
\mathbf{T} (\mathcal{N}, \mathcal{N}) =& \frac{1}{2 r^4 \sin^4 \theta}
\end{align}

The energy of the resulting wave equation is 
\begin{align}
E(t) =& \int_{\Sigma_t} \mathbf{T} (\mathcal{N}, \mathcal{N})\sqrt{q} \, dr d\theta \notag\\
& \to \infty.
\end{align}

\noindent The fact that the simplest solution arising out of the axisymmetric dimensional reduction of $\mathbb{R}^{3+1}$ exhibits such a rampant asymptotic behaviour, that neither correlates with its parent $\mathbb{R}^{3+1}$ nor with $\mathbb{R}^{2+1}$ when it is supposed to help in our construction as a ground-state solution, shall have important implications in what is to follow. Let us now consider Schwarzschild.

\subsection*{Dimensional Reduction of Schwarzschild}

As an example also consider the 
Schwarzschild metric 
\begin{align}
\bar{g}_{sh} = -\left( 1- \frac{2m}{r} \right) dt^2 + \left( 1-\frac{2m}{r} \right)^{-1} dr^2 + r^2 d\omega^2_{\mathbb{S}^2}.
\end{align}
For simplicity suppose $f(r) = \left( 1-\frac{2m}{r} \right),$  then the spatial part is 
\begin{align}
\bar{q}_{sh} = f^{-1} dr^2 + r^2 d\omega^2_{\mathbb{S}^2}
\end{align}

Subsequently,  the reduced $q_{sh}$ can be expressed as
\begin{align}
\bar{q}_{sh} =r^2 \sin^2 \theta \left( f^{-1} dr^2 + r^2 d\theta^2  \right).
\end{align}
The Schwarzschild spacetime is also an example of a $\ptl_\phi$ hypersurface orthogonal spacetime, therefore the
reduced equations are: 
\begin{align}
\mathbf{E}_{\mu \nu} =& \mathbf{T}_{\mu \nu} \quad \textnormal{on} \quad (M, g) \notag \\
\square_g u =&0 \quad \quad \,\, \textnormal{on} \quad (M, g)
\end{align}
where 
\[ g = r^2 \sin^2 \theta ( -f dt^2 + f^{-1} dr^2 + r^2 d\theta^2) \]
and
\[\mathbf{T}_{\mu \nu} = \grad_\mu u \grad_\nu u - \halb g_{\mu \nu} \grad^\sigma u \grad_\sigma u,\]
 $u= \log r \sin\theta,$ like before. The parameter $m$ coincides with $m_{ADM}$ defined earlier.  The curvature of $(\Sigma_{sh}, q_{sh})$ also blows up at infinity and its geometric-mass diverges. Again consider $\mathbf{T} (\mathcal{N}, \mathcal{N})$ for $\mathcal{N} = (r\sin\theta\sqrt{f} )^{-1} \ptl_t $ for $r$ sufficiently large. 
We have, 
\begin{align}
\grad^\sigma u \grad_\sigma u =& g^{rr} \ptl_r u \ptl_u r + g^{\theta \theta} \ptl_\theta u \ptl_\theta u \notag \\
=& f \frac{1}{r^4 \sin^2\theta} + \frac{\cos^2 \theta}{r^4 \sin^4 \theta} \\
=& \frac{1}{r^4 \sin^2 \theta} \left( f + \cot^2 \theta \right) \\
=& \frac{1}{r^4 \sin^2 \theta}  \left( \csc^2 \theta - \frac{2m}{r}  \right)
\end{align}
Subsequently, 
\begin{align}
E(t) =& \int_{\Sigma'} \mathbf{T}(\mathcal{N}, \mathcal{N}) \, \bar{\mu}_q \to \infty.  
\end{align}

\subsection*{Dimensional Reduction of Kerr}

In Boyer-Lindquist coordinates Kerr solution can be represented as

\begin{align}
\bar{g}_{k} = - \left(1-\frac{2Mr}{A} \right)dt^2 - \frac{4Mra \sin^2 \theta}{A} dt d\phi + \frac{B \sin^2 \theta}{A} d\phi^2 + \frac{A}{C} dr^2 + A\,d \theta^2. 
\end{align}

where 

\begin{align}
A = r^2 + a^2 \cos^2 \theta,\, B = \left(r^2 + a ^2 \right)^2 -  C a^2 \sin^2 \theta,\, C = r^2 -2Mr + a^2
\end{align}

As $\ptl_\phi$ is a rotational Killing vector of the Kerr metric, it also admits the $\ptl_\phi$ dimensional reduction. In addition, it also admits stationarity in the asymptotic region, so in this region the time dependence of the wave maps equation drops out. In tune with our discussion, the wave map arising out of $\ptl_\phi$ dimensional reduction of Kerr has
infinite energy
%\begin{align}
%\mathfrak{E}(U) = \int_{\Sigma'} q^{ab} h^{ij} \ptl_a U^i \ptl_b U^k \bar{\mu}_{q}
%\end{align}
\begin{align}
\mathfrak{E}(U) = \int_{\Sigma'} \mathbf{e} \, \bar{\mu}_{q}
\end{align}
This behaviour is irrespective of the issue at the axis and is unaffected by the ergo-region. Indeed, one can prove in general that the AV-mass of $(\Sigma, q)$ obtained by the aforementioned dimensional reduction from $ (\bar{\Sigma}, \bar{q})$ shall always be divergent even though $m_{ADM}$ of $ (\bar{\Sigma}, \bar{q})$  converges. 

Starting from Minkowski to the Kerr family, we have generated  infinite energy (divergent) 2+1 systems. This is already an awry situation but let us analyze this further.   

Firstly, we would like to point out the reason: by construction, the norm of the rotational Killing vector blows up at spatial infinity whereas there is no such restriction on the norm of the translational isometry. Indeed, a calculation shows that the equivalent dimensional reduction of $\mathbb{R}^{3+1}$ using a translational Killing vector does not have any of these issues. 

From a PDE perspective, wave maps are considered to be tractable structures as they are natural geometric generalizations of harmonic maps and linear wave equations. In addition, wave maps have better null structure than Einstein's equations. Since the Kerr family admits the ansatz \eqref{axi-ansatz}, in principle the Kerr family is also a `solution' to the 2+1 Einstein wave map system. Therefore, in theory one can consider the dynamical perturbations within the 2+1 Einstein-wave map class. However, in view of the previously illustrated divergences, this consideration becomes quite subtle. The details are explained below:
\begin{enumerate}

\item When faced with such divergences, a natural resolution that comes to mind is to introduce weights to induce the desired 2+1 asymptotic behaviour. Importantly, we would like to emphasize that any geometrically consistent way of inducing the desired asymptotics is bound to disturb the wave map structure, thus undermining the whole construction centered around their elegance. 

\iffalse
\item Furthermore, even if there is a geometrically consistant way of introducing better behaved (corrected) energies, the global existence theory could become very complicated in these norms
% as these corrections could potentially tamper with the natural scales of the Einstein-wave map system
and certainly would not be comparable to the large-data result  in \cite{AGS}.
\fi
\item In another marked distinction with the equivariant Einstein-wave map system, the divergent Einstein-wave map energies mean that the 2+1 wave map stress-energy fluxes are not the honest carriers of the true stress-energy fluxes propagating in the parent 3+1 axisymmetric system.

\item An issue that is somewhat independent of the aforementioned issues is the general issue of proving the decay of waves in 2+1 dimensions. This can be traced back to the failure of classical Huygens principle for waves in $\mathbb{R}^{2+1}$ and the associated weaker decay. In view of the fact that the major open problems on stability of blackholes are directly related to optimal decay rates, this issue also assumes significance. 
\end{enumerate}
In a forthcoming work we shall attempt to overcome some of these issues using a Hamiltonian approach for the \emph{linearized} framework of the problem. 
%Consistant with these observations is the fact that the dimensional reduction approach has almost never been adopted for dynamical axisymmetric spacetimes in the past (see in contrast \cite{})

\section*{Ernst Equations}
Separately, we would like to emphasize that the dynamical axisymmetric dimensional reduction has to be strictly contrasted with a related system of equations that admit the Schwarzschild and Kerr spacetimes as solutions, but arrived via stationarity asumption and subsequent imposition of axisymmetry. 

\noindent Suppose $(\bar{M}, \bar{g})$ is a 3+1 dimensional globally hyperbolic Lorentzian spacetime. $(M, g)$ is called a stationary spacetime if there exits a timelike Killing vector $T$ outside a compact set of the Cauchy hypersurfaces $\olin{\Sigma}_S, \bar{q}_S$. In our construction, in order to consider the physically relevant cases we would of course omit the case with orbiting or closed integral curves of $T$. Then consider the ansatz \footnote{variants of this ansatz are called Weyl coordinates and Papapetrou coordinates}

\begin{align}
\bar{g} = \bar{q}_S + e^{2u_S} \, \mathbb{T}^2 
\end{align} 
where $\mathbb{T} \fdg = dt + (A_S)_a dx^a,\, a = 1, 2, 3.$

The vacuum Einstein's equations in this case 
\begin{align}
\bar{R}_{\mu \nu} =0 \quad \textnormal{on} \quad (M, g)
\end{align} 
can also be reduced in essentially a similar way,  but now by using the time-like Killing vector from the Lorentzian 4-manifold to a Riemannian 3-manifold $(\bar{\Sigma}, \bar{q})$ (see Geroch's projection formalism \cite{geroch-1971, geroch-1972}) outside the ergoregion and a subsequent conformal transformation
\[ \tilde{q}_S \fdg = e^{2 u_S} \,\, \bar{q}_S \]
As before, we shall denote $\tilde{q}_S$ by $\bar{q}_S$ itself for simplicity.

Now, given the reduced stationary Einstein's equations on $(\bar{\Sigma}_S, \bar{q}_S)$, we can impose axisymmetry: Suppose there exists a rotational Killing vector $\ptl_{\phi},$ with closed orbits and fixed points. Then consider 
$(\Sigma_S, q_S) \fdg = (\olin{\Sigma}_S, \bar{q}_S) \setminus SO(2)$, we shall arrive at the following system of equations 

\begin{align}
\Delta_{q_S} \Psi^i + \leftexp{(h)}{\Gamma}^i_{jk}\, q_{S}^{a b} \ptl_a \Psi^j \, \ptl_b \Psi^k =0,
\end{align}

$\Psi \fdg (\Sigma_S, q_S) \to (\mathbb{H}^2, h),$ where $(\mathbb{H}^2, h)$ is again the hyperbolic 2-plane, $\Gamma$'s are its Christoffel symbols and $\Delta_q \fdg$ is the Laplace-Beltrami operator 
corresponding the covariant derivative of $(\Sigma, q).$ $\Psi = (u_S, v_S)$ where  $u_S$ is related to the lapse and the twist $v_S$ is constructed from the shift vector.

Aside from the expected issue at the axis, this harmonic map is a well-defined PDE object and one can easily specify conditions so that it is of finite energy. Historically, these equations are called Ernst equations \cite{ernst-1, ernst-2,dkramer-1987,liv-rev-MH} \footnote{sometimes more generally without axisymmetry}. 

If one recalls the Kerr metric

\begin{align}
\bar{g}_{k} = - \left(1-\frac{2M r}{A} \right)dt^2 - \frac{4Mra \sin^2 \theta}{A} dt d\phi + \frac{B \sin^2 \theta}{A} d\phi^2 + \frac{A}{C} dr^2 + A \,d \theta^2, \notag\\
A = r^2 + a^2 \cos^2 \theta,\, B = \left(r^2 + a ^2 \right)^2 -  C a^2 \sin^2 \theta,\, C = r^2 -2Mr + a^2,
\end{align}

one can see that the Kerr spacetime also admits the time like Killing vector the stationary reduction in the spatial asymptotic region.

In other words, even though the Kerr solution fits the ans\"atze for both the dimensional reductions, we have arrived at maps with very different properties. This means that the dimensional reduction procedures are not commutable in general \footnote{Indeed, this notion is fundamental in the construction of the Geroch group}.% irrespective
%of the commutator of the Killing vectors. 

These stationary equations, apart from their role during the discovery of Kerr spacetime, have been instrumental in several mathematical and physical results of much value and delight throughout the history. In addition, one can establish a nice relation between $m_{AV}$ of $(\Sigma_S, q_S)$ and $m_{ADM}$ of $(\olin{\Sigma}_S, \bar{q}_S).$ Therefore, it is perhaps advisable to clearly distinguish the Ernst equations with \eqref{axi-wm}. %If these intricacies related to dimensional reduction are not highlighted, it could lead to the misconception that it is natural to study 2+1 dynamical wave map (arising from dynamical axisymmetric dimensional reduction, which results in  infinite energy) perturbations of the Ernst equations of the Kerr family. \\
\subsection*{Acknowledgements} The author gratefully acknowledges the conversations with Abhay Ashtekar.
 \bibliography{main.bib}
\bibliographystyle{plain}
\emph{[The author wrote this note in the interest of science and to be of service to the mathematical general relativity community]}
\end{document}